\documentclass{article}
\usepackage[pagewise]{lineno}
\usepackage[pdftex, dvips]{graphicx} 
\graphicspath{{images/}{../images/}}
\usepackage{subfig}
\DeclareGraphicsExtensions{ .eps, .png}
\usepackage{float}
\usepackage[english]{babel}
\usepackage{latexsym,amsmath,amssymb,amsfonts,cancel}
\usepackage{bbm}
\usepackage{hyperref}
\usepackage{fancyhdr}
\usepackage[margin=2.5cm,top=3cm]{geometry}

\usepackage{cite}
\usepackage{csquotes}
\usepackage{multirow}

\newtheorem{theorem}{Theorem}

\newtheorem{corollary}{Corollary}
\newtheorem{remark}{Remark}

\newtheorem{definition}{Definition}

\def\R{{\rm I}\hskip -0.85 mm{\rm R}}
\def\N{{\rm I}\hskip -0.85 mm{\rm N}}

\begin{document}

\title{ \huge 
\textbf{Some questions about the regularity and the uniqueness of solutions of parabolic partial differential equations}}

\author{ Inmaculada Gayte-Delgado\thanks{Dpto. de Ecuaciones Diferenciales y An\'alisis Num\'erico, Universidad de Sevilla, Seville, Spain. E-mail: gayte@us.es}; Irene Mar\'{\i}n-Gayte\thanks{Departamento M\'etodos Cuantitativos, Universidad Loyola Andaluc\'ia, Campus Sevilla, Sevilla, Spain.  E-mail:imgayte@uloyola.es}
 } 

\date{}

\maketitle

 \begin{abstract}
This work obtains a fixed-point equation for the solution of linear parabolic partial differential problems based on solutions to heat problems. This is a pointwise equality, so we have required non-standard techniques that involve the study of the sign of certain solutions to linear parabolic problems. This fixed-point equation implies regularity properties of solutions to parabolic problems, not necessarily linear, and this allows us to prove the uniqueness of the solution in three dimensions for the Navier-Stokes problem.
 \end{abstract}
 
 \vspace*{0,5in}
 
\textbf{AMS Classifications:} 35Q30, 49K20, 76D05, 93C95.

\textbf{Keywords:} Initial-boundary value problems for second-order parabolic equations, Navier-Stokes equations, uniqueness problems for PDEs, Smoothness and regularity of solutions to PDEs, Weak solutions to PDEs.

\section{Introduction}

Partial differential equations (PDEs) of parabolic type play a fundamental role in the mathematical modeling of diffusion, heat propagation, and viscous fluid dynamics. These equations are characterized by the presence of a time derivative and a second-order spatial operator, typically the Laplacian, and are known for their smoothing properties and the well-posedness of initial-boundary value problems under suitable assumptions.

Among the most important nonlinear parabolic systems are the incompressible Navier-Stokes equations, which describe the motion of viscous fluid flows. In their classical form, the Navier-Stokes equations couple a momentum equation with an incompressibility constraint, and their analysis poses significant mathematical challenges. 

In two dimensions, the theory is more developed: weak solutions exist and are unique due to the better control over the nonlinear convective term. Key results can be found in the classical works of Lions~\cite{Lions}, Temam~\cite{Teman}, and Prodi~\cite{Prodi}. In this case, there are good estimations of the nonlinear term in the equation which makes the weak solution is smooth enough. 
However, in three dimensions, the uniqueness is assured in a set of smooth solutions, although weaker than in two dimension, but it is not assured the existence of solution in this set.  If it is supposed that the data are small enough or the viscosity is big enough, then the regularity of solutions is guaranteed and so, there is a unique solution in three dimension.

In the study of the three-dimensional Navier?Stokes equations, one of the first approaches consisted in extending the results known in the two-dimensional case to the three-dimensional setting, see \cite{Teman, Lions, HeyRanTur96, Ladyzhenskaya1}. However, all these attempts encountered a major difficulty: the need for more regular solutions, whose existence is not yet guaranteed, (they need a smooth solution in $L^\infty (0,T; L^4(\Omega))$),   or alternatively, the requirement of smoother initial data.

A significant contribution in this direction was made by Caffarelli, Kohn, and Nirenberg, who introduced a local energy inequality as part of their attempt to prove uniqueness of weak solutions in three dimensions, see \cite{Caffarelli}. This inequality includes a term representing the physical energy dissipated over time due to viscous friction, often referred to as viscous energy. It can be interpreted as a local estimate indicating that, at each point in space, the energy lost due to friction is bounded by terms involving only the velocity, pressure, and external forces acting on the fluid.

The local energy inequality was introduced for two main reasons. Firstly, it provides a rigorous formulation of the fact that energy dissipation is a local phenomenon. Secondly, the authors hoped that such an inequality would eventually lead to a uniqueness result for a class of solutions they referred to as suitable  weak solutions. However, this goal remains unachieved to date. Instead, their main result proves that the set of singular points, where the solution may fail to be regular, has zero one-dimensional Hausdorff measure.

An alternative and significantly simpler proof of this regularity result was later given by Lin, see \cite{Lin}.  His argument relies mainly on the estimates developed by Sohr and Wahl in \cite{HS-WW},  and leads to the same conclusion as the original work of Caffarelli, Kohn, and Nirenberg.

Furthermore, \cite{Caffarelli} also includes an existence theorem for suitable weak solutions under a specific approximation framework. This method is based on a time-discretization scheme where each step uses the solution from the previous time level. Related developments can be found in other works that provide new proofs of the Caffarelli-Kohn-Nirenberg theorem. In particular, Guermond  \cite{JLG} employs a different Galerkin-type scheme that reaches the same result, while Da Veiga \cite{DaVeiga} presents an alternative proof based on a fourth-order regularization method.

In this work, we prove the uniqueness of the solution to the Navier-Stokes problem in three dimensions by following the reasoning of \cite{Teman}. To this purpose, we introduce a new analysis of parabolic PDEs. These tools allow us to obtain a new characterization of the solution of a linear parabolic problem for each $t$. In particular, we show that the solutions of linear parabolic problems lie in $L^\infty(0,\,T;L^4(\Omega))$ when the initial data is in $L^4(\Omega)$. Therefore, we can establish the uniqueness of solutions to the Navier-Stokes equations in three dimensions.

In the first section, we focus on the main result which is based on the proof of uniqueness of the Navier-Stokes problem:
\begin{theorem}
Let be $u:\Omega\times (0,\,T)\to \R$ sastisfying $0<c\le u\le M$, $y_0\in L^2(\Omega)$ with $y_0\ge 0$ non identically zero, and $y$ the solution of 
$$\left\{
\begin{array}{l}
y_t-\Delta y=u\hbox{ in }\Omega\times (0,\,T)\\
\noalign{\smallskip}
y=0 \hbox{ on }\partial \Omega\times (0,\,T)\\
\noalign{\smallskip}
y(0)=y_0\hbox{ in }\Omega.
\end{array}
\right.
$$
Then, the following equality is stated
$$y(t)=\displaystyle{\frac{\vert y(t)\vert}{\vert \varphi(t)\vert}}\varphi(t)\ \ \ \forall t\in [0,\,T],$$
being $\Omega \subset \R^N$ a bounded open set, with boundary in $C^{0,1}$, $\vert\cdot\vert$ the norm in $L^2(\Omega)$ and $\varphi$ the solution of the heat problem
$$\left\{
\begin{array}{l}
\varphi_t-\Delta \varphi=0\hbox{ in }\Omega\times (0,\,T)\\
\noalign{\smallskip}
\varphi=0 \hbox{ on }\partial \Omega\times (0,\,T)\\
\noalign{\smallskip}
\varphi(0)=y_0\hbox{ in }\Omega.
\end{array}
\right.
$$
\end{theorem}
This result establishes that $y(t)$ has the same spatial regularity as $\varphi(t)$ and it is continuous in $[0,\,T]$. The idea of the proof of this result comes from \cite{IreneGayte} where it is proved the existence of controls which are partially distributed in linear parabolic equations with general diffusion coefficients. First of all, it is shown the existence of functions $v$ such that, if $\Psi_v$ is the solution of
$$\left\{
\begin{array}{l}
(\Psi_v)_t-\Delta \Psi_v=v\hbox{ in }\Omega\times (0,\,T)\\
\noalign{\smallskip}
\Psi_v=0 \hbox{ on }\partial \Omega\times (0,\,T)\\
\noalign{\smallskip}
\Psi(0)=y_0\hbox{ in }\Omega,
\end{array}
\right.
$$
then 
$$\vert y(T)\vert\Psi_v(T)-y(T)\vert \Psi_v(T)\vert$$ 
has a constant sign in $\Omega$. This implies immediately that
$$\vert y(T)\vert\Psi_v(T)-y(T)\vert \Psi_v(T)\vert=0\hbox{ in }\Omega,$$ when the data are non negative.\newline
We prove that this function $v$ can be chosen zero and so,
$$\vert y(T)\vert\varphi(T)-y(T)\vert \varphi(T)\vert=0\hbox{ in }\Omega.$$

In Section 3 we want to extend the previous result to the case $u\in L^2(0,\,T;H^{-1}(\Omega))$ and non constraint on the sign of $u$ and $y_0$. In this case, the following theorem is hold:
\begin{theorem}
Let be $u\in L^2(0,\,T;H^{-1}(\Omega))$, $y_0\in L^2(\Omega)$. Then, there exist $\lambda_1,\, \lambda_2\in C([0,\,T];(1,+\infty))$ such that
$$y(t)=\lambda_1(t)\varphi_1(t)-\lambda_2(t)\varphi_2(t),$$
with $\varphi_1$ and $\varphi_2$ two solutions of the heat equation with homogeneous Dirichlet boundary conditions.
\end{theorem}
By this theorem we obtain that, in a linear parabolic problem, when the right-hand side is in $L^2(0,\,T;H^{-1}(\Omega))$ and the initial data is in $L^4(\Omega)$, then the solution is in $L^\infty(0,\,T;L^4(\Omega))$.\newline
The last section addresses the uniqueness of the solution of the Navier-Stokes problem in three dimension. When the initial data of the Navier-Stokes problem is in $H\cap L^4(\Omega)^N$, being $H$ the usual space in Navier-Stokes theory,
$$H=\{v\in L^2(\Omega)^N:\ \  \nabla \cdot v=0\hbox{ in }\Omega,\ v\cdot \vec{n}=0\hbox{ on }\partial \Omega\},$$
and the right-hand side in $L^2(0,\,T;H^{-1}(\Omega)^N)$, the uniqueness is obtained writing the solution of Navier-Stokes problem as the sum of a smooth function and a solution of a linear parabolic problem. Since it has been proved the integrability property for the solutions of linear parabolic problems, we have that the solution of Navier-Stokes problem is in $L^\infty(0,\,T;L^4(\Omega)^N)$. Next we will relax the hypotheses: first of all, if the initial data is in $H$, using a discretization through the eigenfunctions of the Stokes operator, we have a sequence of Navier-Stokes problems with initial data in $H\cap L^4(\Omega)^N$. These problems have a unique solution, and the whole sequence of solutions is shown that converges to a solution of Navier-Stokes with initial data in $H$. So, any solution is reached by the limit of the sequence of approximated solutions, and it is the whole sequence that converges. This implies the uniqueness of solution when the initial data is in $H$ and the right-hand side is in $L^2(0,\,T;H^{-1}(\Omega)^N)$.\newline
The last step consists of proving the uniqueness when the right-hand side is in $L^2(0,\,T;V')$, being the space $V'$ the dual of $V$:
$$V=\{ v\in H^1_0(\Omega)^N:\ \ \ \nabla\cdot v=0\hbox{ in }\Omega\}.$$

The paper is divided in three sections: in Section 2 we obtain a relation between the solution of the heat problem and the solution of a linear parabolic problem with non negative data. In Section 3 it is studied the case of a linear parabolic problem with any data and we obtain integrability properties of the solution. And, in the last section, we approach the study of the uniqueness of Navier-Stokes solution in three dimension.

\section{A linear parabolic problem with non negative data.}
In this section we focus  on  the study of a linear parabolic problem. The problem is
\begin{equation}\label{pby}
\left\{
\begin{array}{l}
y_t-\Delta y=u\hbox{ in }\Omega\times (0,\,T)\\
\noalign{\smallskip}
y\vert_{\Sigma}=0\\
\noalign{\smallskip}
y(0)=y_0\hbox{ in }\Omega,
\end{array}
\right.
\end{equation}
where $\Omega\subset \R^N$ is a bounded open set whose boundary is locally Lipschitz, $\Sigma=\partial \Omega\times (0,\,T)$, $u\in L^2((0,\,T)\times \Omega)$ and $y_0\in L^2(\Omega)$.\newline
We call $\varphi$ the solution of the heat problem:
\begin{equation}\label{pbvarphi}
\left\{
\begin{array}{l}
\varphi_t-\Delta \varphi=0\hbox{ in }\Omega\times (0,\,T)\\
\noalign{\smallskip}
\varphi\vert_{\Sigma}=0\\
\noalign{\smallskip}
\varphi(0)=y_0\hbox{ in }\Omega,
\end{array}
\right.
\end{equation}
and $W(0,T)$ the space
$$W(0,T)=\{ w\in L^2(0,T;H^1_0(\Omega)):\ \ w_t\in L^2(0,T;H^{-1}(\Omega))\}.$$
We denote $(\cdot,\cdot)$ the scalar product in $L^2(\Omega)$, $<\cdot,\cdot>$ the scalar product in $L^2(\Omega\times (0,\,T))$ and $\vert \cdot\vert$ the norm in $L^2(\Omega)$.

We begin with a result which can be considered a maximum principle at the final time and it is proved in \cite{IreneGayte}. We show it here for completeness.

\begin{theorem}\label{Ppiomaximo}
Let be $\beta\in C^1([0,T])$ verifying
\begin{equation}\label{propiedadesbeta}
\begin{array}{l}
\beta>0\\
\noalign{\smallskip}
\max_{[0,\,T]}\beta=\beta(T),
\end{array}
\end{equation}
$w\in W(0,T)$ such that
\begin{equation}\label{propiedadesw}
\left\{
\begin{array}{l}
w_t-\Delta w\ge 0 \hbox{ in }\Omega\times (0,\,T)\\
\noalign{\smallskip}
w\vert_{\Sigma}=0 \\
\noalign{\smallskip}
w(0)\ge 0 \hbox{ in }\Omega.
\end{array}
\right.
\end{equation}
Then, the solution $z$ of the problem 
\begin{equation}\label{pbz}
\left\{
\begin{array}{l}
z_t-\Delta z=-\beta'w \hbox{ in }\Omega\times (0,\,T)\\
\noalign{\smallskip}
z\vert_{\Sigma}=0 \\
\noalign{\smallskip}
z(0)=z_0\hbox{ in }\Omega
\end{array}
\right.
\end{equation}
with $z_0\le 0$ in $\Omega$ verifies 
$$z(T)\le 0.$$
\end{theorem}
\begin{remark}
If $\beta$ is increasing, the statement of Theorem \ref{Ppiomaximo} follows from the classical comparison theorem since in this case $-\beta'w\le 0$. The main contribution of the Theorem \ref{Ppiomaximo} is that $\beta$ is not increasing necessarily in $(0,\,T)$.
\end{remark}

\textbf{Proof:}
Let  $\tilde{\Omega}\subset \Omega$ be any open set. Consider the backward heat problem
$$
\left\{
\begin{array}{l}
-\eta_t-\Delta \eta=0\hbox{ in }\Omega\times (0,\,T)\\
\noalign{\smallskip}
\eta\vert_{\Sigma}=0\\
\noalign{\smallskip}
\eta(T)=\displaystyle{1_{\tilde{\Omega}}}\hbox{ in }\Omega,
\end{array}
\right.
$$
being $\displaystyle{1_{\tilde{\Omega}}}$ the characteristic function in $\tilde{\Omega}$, i.e.
$$\displaystyle{1_{\tilde{\Omega}}}=
\left\{
\begin{array}{l}
1\hbox{ if }x\in \tilde{\Omega}\\
\noalign{\smallskip}
0\hbox{ if }x\in \Omega\setminus\tilde{\Omega}.
\end{array}
\right.
$$
Multiplying the partial differential equation of (\ref{pbz}) by $\eta$ and integrating by parts, we obtain 
\begin{equation}\label{eqsign1}
(z(T),1_{\tilde{\Omega}})=(z_0,\eta(0))+<-\beta'w,\eta>.
\end{equation}
Let
$$F(s)=s-\beta(T).$$
Then
$$<-\beta'w,\eta>=<-\frac{d\,F(\beta(t))}{dt}w,\eta>=\int_0^T-\frac{d}{dt}F(\beta(t))(w,\eta).$$
Integrating by parts this last integral we get
$$<-\beta'w,\eta>=\int_0^TF(\beta(t))\frac{d}{dt}(w,\eta)\, -\, F(\beta(T))(w(T),\eta(T))+F(\beta(0))(w(0),\eta(0)).$$
Since $\beta$ attains its maximum in $t=T$, we obtain 
$$F(\beta(t))\le 0\ \ \ \forall t\in [0,\,T],$$
and besides,
$$F(\beta(T))=0,$$
so
$$<-\beta'w,\eta>\le\int_0^TF(\beta(t))\frac{d}{dt}(w,\eta).$$
Next, we study 
$$\frac{d}{dt}(w,\eta).$$
$$\frac{d}{dt}(w,\eta)=(w_t,\eta)+(w,\eta_t)=(w_t,\eta)-(w,\Delta\eta)=$$
$$=(w_t-\Delta w,\eta)\ge 0,$$
and, since $F(\beta(t))\le 0$ in $[0,\,T]$
$$<-\beta'w,\eta>\le 0.$$
Using this last inequality, (\ref{eqsign1}) and $z_0\le 0$ we obtain
$$\displaystyle{\left(z(T),1_{\tilde{\Omega}}\right)}\le 0,$$
i.e., $z(T)\le 0$ in $\Omega$ since $\tilde{\Omega}$ is any arbitrary open set in $\Omega$.
\begin{flushright}
$\square$
\end{flushright}

Next, the main result of this section is stated:
\begin{theorem}\label{multiplodevarphi}
Let be $y_0\ge 0$ non identically zero, $u\in L^\infty(\Omega\times (0,\, T))$, $0<c\le u\le M$ in $\Omega\times (0,\,T)$.\newline
Then, it is verified that
$$\displaystyle{\frac{y(T)}{\vert y(T)\vert}}=\displaystyle{\frac{\varphi(T)}{\vert \varphi(T)\vert}},$$
being $y$ the solution of (\ref{pby}) and $\varphi$ the solution of (\ref{pbvarphi}).
\end{theorem}
\textbf{Proof:}
\vskip0.2cm
\noindent
{\bf Step 1:}\newline 
We prove that for any $\varepsilon\in (0,\,1)$, for any open set $\omega\subsetneq \Omega$ and for any interval $I\subsetneq (0,\,T)$, it is verified that there exists $v_\varepsilon^\ast$ such that
$$v_\varepsilon^\ast=\varepsilon u\hbox{ in }\omega\times I,$$
$$0\le v_\varepsilon^\ast \le u \hbox{ in }\Omega \times (0,\,T)$$
and
$$\vert y(T)\vert \Psi_{v_\varepsilon^\ast}(T)-y(T)\vert \Psi_{v_\varepsilon^\ast}(T)\vert=0\hbox { in }\Omega,$$ 
where $\Psi_{v_\varepsilon^\ast}$ is the solution of
$$
\left\{
\begin{array}{l}
(\Psi_{v_\varepsilon^\ast})_t-\Delta \Psi_{v_\varepsilon^\ast}=v_\varepsilon^\ast\hbox{ in }\Omega\times (0,\,T)\\
\noalign{\smallskip}
(\Psi_{v_\varepsilon^\ast})\vert_{\Sigma}=0\\
\noalign{\smallskip}
\Psi_{v_\varepsilon^\ast}(0)=y_0 \hbox{ in }\Omega.
\end{array}
\right.
$$
For that, we shall build a sequence $\{v_k\}_k$ verifying
$$\vert\Psi_{v_{k-1}}(T)\vert y(T)-\vert y(T)\vert \Psi_{v_k}(T)\le 0\hbox{ in }\Omega,$$
where $\Psi_{v_k}$ is the solution of the problem
\begin{equation}\label{pbPsi}
\left\{
\begin{array}{l}
(\Psi_{v_k})_t-\Delta \Psi_{v_k}=v_k\hbox{ in }\Omega\times (0,\,T)\\
\noalign{\smallskip}
(\Psi_{v_k})\vert_{\Sigma}=0\\
\noalign{\smallskip}
\Psi_{v_k}(0)=y_0 \hbox{ in }\Omega.
\end{array}
\right.
\end{equation}
Let $v_1\in L^\infty(\Omega\times (0,T))$, $0\le v_1\le u$, such that
\begin{equation}\label{v0noesu}
v_1(x,t)=\varepsilon u(x,t)\ \ \ \hbox{for a.e. } (x,t)\in \omega\times I.
\end{equation}
We choose $\beta_1\in C^1([0,\,T])$ and $w_1\in W(0,T)$ satisfying (\ref{propiedadesbeta}) and (\ref{propiedadesw}), respectively (the hypothesis of Theorem \ref{Ppiomaximo}), and besides,
$$w_1\in L^\infty(\Omega\times (0,T)),
$$
\begin{equation}\label{estimacionbeta1derivada}
\beta'_1 \le \frac{(\vert y(T)\vert -\vert \Psi_{v_1}(T)\vert)c}{\sup_{\Omega \times (0,T)} w_1}\ \ \ \hbox{ in } [0,T],
\end{equation}
\begin{equation}\label{estimacionestrictabeta1derivada}
\beta'_1<-\frac{\vert\Psi_{v_1}(T)\vert \sup_{\omega\times I}u}{\inf_{\omega \times I} w_1} \ \ \ \hbox{ in } I.
\end{equation}
The condition (\ref{estimacionbeta1derivada}) is possible because
$$\displaystyle{\frac{(\vert y(T)\vert -\vert \Psi_{v_1}(T)\vert)c}{\sup_{\Omega \times (0,T)} w_1}}>0.$$
The inequality (\ref{estimacionestrictabeta1derivada}) has sense because ${\inf_{\omega \times I} w_1}>0$ since $w_1$ is strictly bigger than the solution of the heat equation with initial data $w_1(0)$. Besides, this inequality requires that $\beta_1$ be decreasing in $I$ and this is possible because $I\subsetneq (0,\,T)$.\newline
Then, it is verified that
\begin{equation}\label{propiedad1}
\beta'_1 w_1\le (\vert y(T)\vert -\vert \Psi_{v_1}(T)\vert)u\ \ \ \hbox{ in } \Omega\times [0,T].
\end{equation}
Effectively, by (\ref{estimacionbeta1derivada}),
$$\beta_1'\sup_{\Omega \times (0,T)} w_1\le (\vert y(T)\vert -\vert \Psi_{v_1}(T)\vert)c,$$
if $\beta_1'(t)>0$, we have that 
$$\beta_1'(t)w_1(x,t)\le \beta_1'(t)\sup_{\Omega \times (0,T)} w_1,$$
and if $\beta_1'(t)\le0$, then (\ref{propiedad1}) is obvious because $(\vert y(T)\vert -\vert \Psi_{v_1}(T)\vert)u>0$.\newline
It is also verified that
\begin{equation}\label{propiedad2}
\beta'_1(t)w_1(x,t)<-\vert\Psi_{v_1}(T)\vert u(x,t)\ \ \ \forall (x,t)\in \omega\times I.
\end{equation}
because, by (\ref{estimacionestrictabeta1derivada}),
$$\beta_1'\inf_{\omega \times I} w_1<-\vert\Psi_{v_1}(T)\vert \sup_{\omega\times I}u\ \ \ \hbox{ in }I,$$
and since $\beta_1'$ is negative in $I$, we obtain
$$\beta_1' w_1<-\vert\Psi_{v_1}(T)\vert \sup_{\omega\times I}u\ \ \ \hbox{ in }\omega\times I.$$
Note that the property (\ref{propiedad1}) would not be possible if $\vert\Psi_{v_1}(T)\vert=\vert y(T)\vert$, but this does not happen because $v_1\le u$, and this inequality is strict in $\omega\times I$.\newline
We define $\tilde{v}_1$ as
\begin{equation}\label{definicionvtilde}
\tilde{v}_1=\frac{\vert\Psi_{v_1}(T)\vert u+\beta_1'w_1}{\vert y(T)\vert}.
\end{equation}
Then
$$\vert \Psi_{v_1}(T)\vert u -\vert y(T)\vert \tilde{v}_1=-\beta_1'w_1.$$
If we call
$$\tilde{z}_1=\vert \Psi_{v_1}(T)\vert y-\vert y(T)\vert \Psi_{\tilde{v}_1},$$
where $\Psi_{\tilde{v}_1}$ is the solution of a problem like (\ref{pbPsi}) with right-hand side $\tilde{v}_1$, we have
$$
\left\{
\begin{array}{l}
\tilde{z_1}_t-\Delta\tilde{z_1}=-\beta'_1w_1\\
\noalign{\smallskip}
\tilde{z_1}\vert_{\Sigma}=0\\
\noalign{\smallskip}
\tilde{z_1}(0)=(\vert \Psi_{v_1}(T)\vert-\vert y(T)\vert)y_0,
\end{array}
\right.
$$
whose initial data is less than or equal to zero. By Theorem \ref{Ppiomaximo} $\tilde{z_1}(T)\le 0$, i.e.
\begin{equation}\label{desig}
\vert \Psi_{v_1}(T)\vert y(T)-\vert y(T)\vert \Psi_{\tilde{v_1}}(T)\le 0.
\end{equation}
On the other hand, as
$$\tilde{v}_1=\frac{\vert\Psi_{v_1}(T)\vert}{\vert y(T)\vert}u+\frac{\beta_1'w_1}{\vert y(T)\vert},$$
applying (\ref{propiedad1}) to the last term, we get
$$\tilde{v_1}\le 
 \frac{\vert\Psi_{v_1}(T)\vert}{\vert y(T)\vert}u+\frac{\vert y(T)\vert -\vert\Psi_{v_1}(T)\vert}{\vert y(T)\vert}u= u.$$
And by (\ref{propiedad2}) we also get
$$\tilde{v}_1(x,t)< \frac{\vert\Psi_{v_1}(T)\vert}{\vert y(T)\vert}u- \frac{\vert\Psi_{v_1}(T)\vert}{\vert y(T)\vert}u=0\hbox{ for a.e. }(x,t)\in \omega\times I.$$
Then, we define 
\begin{equation}\label{definicionv1}
v_2=\max(\tilde{v}_1,\, v_1).
\end{equation}
The properties of $v_2$ are:
\begin{equation}\label{propiedad1dev1}
v_2\ge v_1
\end{equation}
$$v_2\le u,$$
since $v_1\le u$ and $\tilde{v_1}\le u$, and because of $\tilde{v_1}<0$ in $\omega\times I$
\begin{equation}\label{v1igualav0}
v_2(x,t)=v_1(x,t) \ \ \ \hbox{for a.e. } (x,t)\in \omega\times I
\end{equation}
and so, 
$$v_2=\varepsilon u\hbox{ in }\omega\times I.$$
Besides, by the maximum principle
$$\vert \Psi_{v_1}(T)\vert y(T)-\vert y(T)\vert \Psi_{v_2}(T)\le \vert \Psi_{v_1}(T)\vert y(T)-\vert y(T)\vert \Psi_{\tilde{v_1}}(T).$$
By (\ref{desig}),
$$\vert \Psi_{v_1}(T)\vert y(T)-\vert y(T)\vert \Psi_{v_2}(T)\le 0\hbox{ in }\Omega.$$
Then, we have
\begin{equation}\label{propiedadesdev2}
\left\{
\begin{array}{l}
0\le v_1\le v_2\le u\\
\noalign{\smallskip}
v_2=v_1=\varepsilon u\hbox{ in }\omega \times I\\
\noalign{\smallskip}
\vert \Psi_{v_1}(T)\vert y(T)-\vert y(T)\vert \Psi_{v_2}(T)\le 0\hbox{ in }\Omega.
\end{array}
\right.
\end{equation}
Repeating the reasoning, now with $v_2$ in the place of $v_1$, choosing $\beta_2$ and $w_2$, and so on, we get a sequence of functions  $\{ v_k\}_{k\ge 1}$ such that
\begin{equation}\label{increasing}
0\le v_k\le v_{k+1}\le u
\end{equation}
\begin{equation}\label{desigualdadnegativa}
\vert\Psi_{v_{k-1}}(T)\vert y(T) -\vert y(T)\vert \Psi_{v_k}(T)\le 0\hbox{ in }\Omega.
\end{equation}
\begin{equation}\label{vkesv0}
v_k=v_{k-1}=\varepsilon u \hbox{ in } \omega\times I.
\end{equation}
Then, there exists the limit of $\{v_k\}_k$, for almost $(x,t)\in \Omega\times (0,\,T)$, denoted by $v_{\varepsilon}^\ast$, and it verifies
$$0\le v^\ast_\varepsilon\le u$$
$$v^\ast_\varepsilon=\varepsilon u\hbox{ in }\omega\times I$$
$$\vert \Psi_{v_\varepsilon^\ast}(T)\vert y(T)-\vert y(T)\vert \Psi_{v_\varepsilon^\ast}(T)\le 0\hbox{ in }\Omega.$$
Since $\vert \Psi_{v_\varepsilon^\ast}(T)\vert y(T)+\vert y(T)\vert \Psi_{v_\varepsilon^\ast}(T)>0$ in $\Omega$ and
$$\displaystyle{\left(\vert \Psi_{v_\varepsilon^\ast}(T)\vert y(T)+\vert y(T)\vert \Psi_{v_\varepsilon^\ast}(T),\vert \Psi_{v_\varepsilon^\ast}(T)\vert y(T)-\vert y(T)\vert \Psi_{v_\varepsilon^\ast}(T)\right)}=0$$
we obtain that
$$\vert \Psi_{v_\varepsilon^\ast}(T)\vert y(T)-\vert y(T)\vert \Psi_{v_\varepsilon^\ast}(T)=0\hbox{ in }\Omega,$$
so
\begin{equation}\label{vepsilonestrella}
\displaystyle{\frac{y(T)}{\vert y(T)\vert}}=\displaystyle{\frac{\Psi_{v_\varepsilon^\ast}(T)}{\vert \Psi_{v_\varepsilon^\ast}(T)\vert}}.
\end{equation}
{\bf Step 2:}\newline
Now we have a sequence $\{v_\varepsilon^\ast\}_{\varepsilon}$  bounded in $L^2(\Omega\times (0,\,T))$. For a subsequence, still denoted by $\varepsilon$, there exists the weak limit in $L^2(\Omega\times (0,\,T))$, $v^\ast$. \newline
We know that
$$v_\varepsilon^\ast=\varepsilon u \hbox{ in }\omega\times I,$$
so
$$v^\ast=0\hbox{ in } \omega\times I.$$
Besides, by the weak convergence of the sequence $\{v_\varepsilon^\ast\}_{\varepsilon}$ to $v^\ast$, the sequence of the states $\{\Psi_{v_\varepsilon^\ast}\}$ converges to $\Psi_{v^\ast}$ in $C([0,\,T];L^2(\Omega))$. Indeed, there exists a subsequence of convex linear combinations of $v_{\varepsilon_n}^\ast$ which converges strongly to $v^\ast$ in $L^2(\Omega\times (0,\,T))$,
$$\sum_{i\in I_k}\lambda_iv_{\varepsilon_{n_i}}^\ast \to v^\ast \hbox{ in }L^2(\Omega\times (0,\,T))\hbox{ when }k\hbox{ tends to }\infty,$$
being $I_k$ a finite set of natural numbers, $\lambda_i\ge 0$ and $\sum_{i\in I_k}\lambda_i=1$.\newline
By (\ref{vepsilonestrella})
$$\vert \Psi_{v^\ast_{\varepsilon_{n_i}}}(T)\vert y(T)-\Psi_{v^\ast_{\varepsilon_{n_i}}}(T)\vert y(T)\vert=0$$
and this implies that
$$\sum_{i\in I_k}\lambda_i\vert \Psi_{v^\ast_{\varepsilon_{n_i}}}(T)\vert y(T)-\sum_{i\in I_k}\lambda_i\Psi_{v^\ast_{\varepsilon_{n_i}}}(T)\vert y(T)\vert =0.$$
The linearity of the problem provides the equality
$$\sum_{i\in I_k}\lambda_i\Psi_{v^\ast_{\varepsilon_{n_i}}}(T)=\Psi_{\sum_{i\in I_k}\lambda_iv^\ast_{\varepsilon_{n_i}}}(T),$$
and the convexity of the norm implies that
$$\displaystyle{\vert \sum_{i\in I_k}\lambda_i\Psi_{v^\ast_{\varepsilon_{n_i}}}(T)\vert}\le \sum_{i\in i_k}\lambda_i\vert \Psi_{v^\ast_{\varepsilon_{n_i}}}(T)\vert,$$
so,
$$\displaystyle{\left\vert \Psi_{\sum_{i\in I_k}\lambda_iv^\ast_{\varepsilon_{n_i}}}(T)\right\vert} y(T)-\displaystyle{\Psi_{\sum_{i\in I_k}\lambda_iv^\ast_{\varepsilon_{n_i}}}(T)}\vert y(T)\vert \le 0.$$
Passing to the limit when $k$ tends to infinite we get
$$\vert \Psi_{v^\ast}(T)\vert y(T)-\Psi_{v^\ast}(T)\vert y(T)\vert \le 0\hbox{ in } \Omega.$$
And repeating the reasoning at the bottom of the Step 1 we obtain that
$$\vert \Psi_{v^\ast}(T)\vert y(T)-\Psi_{v^\ast}(T)\vert y(T)\vert =0 \hbox{ in }\Omega.$$
We have just proved that there exists $v^\ast\in L^2(\Omega \times (0,\,T))$ such that
$$0\le v^\ast \le u$$
$$v^\ast=0\hbox{ in }\omega\times I$$
and
$$ \displaystyle{\frac{y(T)}{\vert y(T)\vert}}=\displaystyle{\frac{\Psi_{v^\ast}(T)}{\vert \Psi_{v^\ast}(T)\vert}}.$$
{\bf Step 3:}\newline
Let be a countable open sets covering of $\Omega$, $\{\omega_n\}_n$, $\omega_n\subset \omega_{n+1}$, and a countable open intervals covering of $(0,\,T)$, $\{I_n\}_n$, $I_n\subset I_{n+1}$.\newline
By Step 2, for any $n\in \N$ there exists $v_n^\ast\in L^2(\Omega\times (0,\,T))$ verifying
$$0\le v_n^\ast\le u$$
$$v_n^\ast=0\hbox{ in }\omega_n\times I_n$$
and
$$ \displaystyle{\frac{y(T)}{\vert y(T)\vert}}=\displaystyle{\frac{\Psi_{v_n^\ast}(T)}{\vert \Psi_{v_n^\ast}(T)\vert}}.$$
Repeating the reasoning of the Step 2, there exists a subsequence of $\{v_n^\ast\}$, still denoted by $\{n\}$, such that it converges to $v^\ast$ weakly in $L^2(\Omega\times (0,\,T))$. Then,
$$v^\ast=0 \hbox{ in }\omega_n\times I_n\ \ \ \forall n,$$
so,
$$v^\ast=0\hbox{ in }\Omega\times (0,\,T).$$
Besides,
$$\vert \Psi_{v^\ast}(T)\vert y(T)-\Psi_{v^\ast}(T)\vert y(T)\vert=0.$$
Then, the function $\Psi_{v^\ast}$ is the solution of (\ref{pbvarphi}), that is $\Psi_{v^\ast}=\varphi$.
\begin{flushright}
$\square$
\end{flushright}

\begin{corollary}\label{regularidad}
In the conditions of Theorem \ref{multiplodevarphi} it is verified that $y\in C((0,\,T];C^{\infty}(\Omega))$ and
$$\displaystyle{\frac{y(t)}{\vert y(t)\vert}}=\displaystyle{\frac{\varphi(t)}{\vert \varphi(t)\vert}}\ \ \ \forall t\in (0,\,T].$$
\end{corollary}

\textbf{Proof:}
Consider the problem (\ref{pby}) in $(0,\,t)\times \Omega$ and apply Theorem \ref{multiplodevarphi}.
\begin{flushright}
$\square$
\end{flushright}

\begin{remark}
This result provides a fixed point equation for the solution of a linear parabolic problem sa\-tis\-fying the hypothesis of Theorem \ref{multiplodevarphi}. Besides, it states that this solution in each $t$ is proportional to the solution of the heat problem with the same initial data.\newline
It also states that for any solution $y$ of a linear parabolic problem, with right-hand side $u$ in the hypothesis of Theorem \ref{multiplodevarphi} and initial data $y_0\ge 0$, the quotient $y(t)/\vert y(t)\vert$ for each $t$ is the same function in $\Omega$, $\varphi(t)/\vert\varphi(t)\vert$, which only depends on $y_0$.
\end{remark}

\section{Integrability of the solution of a linear parabolic pro\-blem.}
The goal of this section is to prove that the solution of a linear parabolic problem, when the initial data $y_0$ is in $L^4(\Omega)$, is in $L^\infty(0,\,T;L^4(\Omega))$. For it, we shall show that the solution of a linear parabolic problem, without constraints on the sign of $u$ and $y_0$, is written using the solutions of heat equations.

\begin{theorem}\label{cualquiersigno}
Let be $u\in L^\infty(\Omega\times (0,\,T))$, $y_0\in L^2(\Omega)$. Then, there exists two numbers, $\lambda_{T,1}$ and $\lambda_{T,2}$ in $(1,\,\infty)$, such that
$$y(T)=\lambda_{T,1}\varphi_1(T)-\lambda_{T,2}\varphi_2(T),$$
being $y$ the solution of (\ref{pby}) and $\varphi_1$ and $\varphi_2$ two solutions of the heat equation.
\end{theorem}

\textbf{Proof:}
Let be $c>0$, $u_1=u+\vert u\vert+c$, $y_{01}=y_0+\vert y_0\vert$, $u_2=\vert u\vert+c$  and $y_{02}=\vert y_0\vert$.\newline
We call $y_i$, $i=1,\,2$, the solution of the problem
$$\left\{
\begin{array}{l}
(y_i)_t-\Delta y_i=u_i\hbox{ in }\Omega\times (0,\,T)\\
\noalign{\smallskip}
y_i\vert _\Sigma=0\\
\noalign{\smallskip}
y_i(0)=y_{0i}\hbox{ in }\Omega.
\end{array}
\right.
$$
The functions $u_i$ and $y_{0i}$, $i=1,2$, are in the hypothesis of Theorem \ref{multiplodevarphi} so,
$$y_i(T)=\displaystyle{\frac{\vert y_i(T)\vert}{\vert \varphi_i(T)\vert}}\varphi_i(T),$$
being $\varphi_i$, $i=1,\,2$, the solution of the heat problem
$$\left\{
\begin{array}{l}
(\varphi_i)_t-\Delta \varphi_i=0\hbox{ in }\Omega\times (0,\,T)\\
\noalign{\smallskip}
\varphi_i\vert _\Sigma=0\\
\noalign{\smallskip}
\varphi_i(0)=y_{0i}\hbox{ in }\Omega.
\end{array}
\right.
$$
Denoting 
$$\lambda_{Ti}=\displaystyle{\frac{\vert y_i(T)\vert}{\vert \varphi_i(T)\vert}},\ \ \  i=1,\,2,$$
and using that $y=y_1-y_2$ we obtain the result.
\begin{flushright}
$\square$
\end{flushright}

\begin{theorem}\label{regularidadent}
The statement of Theorem \ref{cualquiersigno} is true for any $t\in (0,\,T]$, i.e., for any $t\in (0,\,T]$ there exists $\lambda_1$, $\lambda_2\in C([0,\,T];(1,\,+\infty))$ such that
$$y(t)=\lambda_1(t)\varphi_1(t)-\lambda_2(t)\varphi_2(t)\ \ \  \forall t\in (0,\,T],$$
being $\varphi_1$ and $\varphi_2$ two solutions of the heat equation. 
\end{theorem}

\textbf{Proof:}
It is enough to consider the problems of the proof of Theorem \ref{cualquiersigno} in $\Omega\times (0,\,t)$.
\begin{flushright}
$\square$
\end{flushright}

Next, we prove an integrability result for the solution of the heat equation whose initial data is in $L^4(\Omega)$.
\begin{theorem}\label{integrabilidadsolucioncalor}
Let be $y_0\in L^4(\Omega)$ and $\varphi$ the solution of the heat problem
$$\left\{
\begin{array}{l}
\varphi_t-\Delta \varphi=0 \hbox{ in }\Omega\times (0,\,T)\\
\noalign{\smallskip}
\varphi\vert_{\Sigma}=0\\
\noalign{\smallskip}
\varphi(0)=y_0\hbox{ in }\Omega.
\end{array}
\right.
$$
Then
$$\varphi\in L^{\infty}(0,\,T;L^4(\Omega)).$$
\end{theorem}

\textbf{Proof:}
The function $\varphi^2$ satisfies the following differential equation in the sense of distributions:
\begin{equation}\label{ecuacionvarphi2}
\displaystyle{\frac{d}{dt}}(\varphi^2)-\Delta (\varphi^2)=2\varphi\varphi_t-2\vert \nabla\varphi \vert ^2-2\varphi \Delta \varphi=-2\vert \nabla\varphi \vert ^2.
\end{equation}
We call $\Psi$ the solution of the heat problem
$$\left\{
\begin{array}{l}
\Psi_t-\Delta \Psi=0\hbox{ in }\Omega\times (0,\,T)\\
\noalign{\smallskip}
\Psi\vert_{\Sigma}=0\\
\noalign{\smallskip}
\Psi(0)=y^2_0\hbox{ in }\Omega.
\end{array}
\right.
$$
By the maximum principle
$$\varphi^2(x,t)\le \Psi(x,t),$$
so
$$\varphi^4(x,t)\le \Psi^2(x,t)$$
and we obtain
$$\Vert \varphi(t)\Vert_{L^4(\Omega)}^4\le \vert\Psi(t)\vert^2_{L^2(\Omega)}.$$
Since $y_0\in L^4(\Omega)$, the solution $\Psi\in C([0,\,T];L^2(\Omega))$, so $\varphi\in L^\infty(0,\,T;L^4(\Omega))$.
\begin{flushright}
$\square$
\end{flushright}

\begin{remark}
The right-hand side of the equation (\ref{ecuacionvarphi2}) is in $L_{loc}^\infty(\Omega\times (0,\,T))$, so we can apply the maximum principle.
\end{remark}

\begin{corollary}
If $u\in L^\infty(\Omega\times (0,\,T))$ and $y_0\in L^4(\Omega)$, then the solution $y$ of the linear parabolic problem (\ref{pby}) is in $L^\infty(0,\,T;L^4(\Omega))$.
\end{corollary}

\textbf{Proof:}
The Theorem \ref{regularidadent} and Theorem \ref{integrabilidadsolucioncalor} imply that when $u\in L^\infty(\Omega\times (0,\,T))$ and $y_0\in L^4(\Omega)$ the solution of a linear parabolic problem is in $L^{\infty}(0,\,T;L^4(\Omega))$.
\begin{flushright}
$\square$
\end{flushright}

Following result states that this is also true when $u$ is more general.
\begin{theorem}\label{Hmenos1L4}
Let be $u\in L^2(0,\,T;H^{-1}(\Omega))$, $y_0\in L^4(\Omega)$. Then, the solution of the problem
$$\left\{
\begin{array}{l}
y_t-\Delta y=u\hbox{ in }\Omega\times (0,\,T)\\
\noalign{\smallskip}
y\vert_{\Sigma}=0\\
\noalign{\smallskip}
y(0)=y_0\hbox{ in }\Omega
\end{array}
\right.
$$
is in $L^{\infty}(0,\,T;L^4(\Omega))$.
\end{theorem}
\textbf{Proof:}
Given $u\in L^2(0,\,T;H^{-1}(\Omega))$, there exists a sequence $\{u_n\}_n\subset L^\infty(\Omega\times (0,\,T))$ such that $\{u_n\}_n$ converges to $u$ the topology of $L^2(0,\,T;H^{-1}(\Omega))$.\newline
By Theorem \ref{regularidadent},
$$y_n(t)=\lambda_{1,n}(t)\varphi_1(t)-\lambda_{2,n}(t)\varphi_2(t),$$
being $y_n$ the solution of 
$$\left\{
\begin{array}{l}
(y_n)_t-\Delta y_n=u_n\hbox{ in }\Omega\times (0,\,T)\\
\noalign{\smallskip}
y_n\vert_\Sigma=0\\
\noalign{\smallskip}
y_n(0)=y_0\hbox{ in }\Omega,
\end{array}
\right.
$$
$\varphi_i$, $i=1,\,2$, the solutions of
$$\left\{
\begin{array}{l}
(\varphi_i)_t-\Delta \varphi_i=0\hbox{ in }\Omega\times (0,\,T)\\
\noalign{\smallskip}
\varphi_i\vert_\Sigma=0\\
\noalign{\smallskip}
\varphi_i(0)=y_{0i}\hbox{ in }\Omega,
\end{array}
\right.
$$
where $y_{0i}$, $y_{i,n}$, $i=1,\,2$ are like in the proof of Theorem \ref{regularidadent}; $y_{01}=y_0+\vert y_0\vert$, $y_{02}=\vert y_0\vert$,
$$\lambda_{i,n}(t)=\frac{\vert y_{in}(t)\vert}{\vert\varphi_i(t)\vert}$$
and $y_{i,n}$ the solution of 
$$\left\{
\begin{array}{l}
(y_{in})_t-\Delta y_{in}=u_{in}\hbox{ in }\Omega\times (0,\,T)\\
\noalign{\smallskip}
y_{in}\vert_\Sigma=0\\
\noalign{\smallskip}
y_{in}(0)=y_{0i}\hbox{ in }\Omega,
\end{array}
\right.
$$
with $u_{1n}=u_n+\vert u_n\vert+c$, $u_{2n}=\vert u_n\vert+c$.\newline
By the continuity of the solution with respect to the data, the states $y_n$ converge when the right-hand sides converge in the topology of $L^2(0,\,T;H^{-1}(\Omega))$. We can pass to the limit when $n$ tends to $\infty$ and we obtain the equality in $L^2(\Omega)$:
$$y(t)=\lambda_1(t)\varphi_1(t)-\lambda_2(t)\varphi_2(t)$$
with
$$\lambda_i(t)=\frac{\vert y_i(t)\vert}{\vert\varphi_i(t)\vert},$$
being $y_i$ the solution of
$$\left\{
\begin{array}{l}
(y_i)_t-\Delta y_i=u_i\hbox{ in }\Omega\times(0,\,T)\\
\noalign{\smallskip}
y_i\vert_\Sigma=0\\
\noalign{\smallskip}
y_i(0)=y_{0i}\hbox{ in }\Omega,
\end{array}
\right.
$$
and $u_1=u+\vert u\vert +c$ and $u_2=\vert u\vert +c$.
Then $y\in L^\infty(0,\,T;L^4(\Omega))$.
\begin{flushright}
$\square$
\end{flushright}

In the case of a system, it remains true the above theorem:
\begin{theorem}\label{sistemaHmenos1L4}
Let be $u\in L^2(0,\,T;H^{-1}(\Omega)^N)$, $y_0\in L^4(\Omega)^N$. Then, the solution of the problem
$$\left\{
\begin{array}{l}
y_t-\Delta y=u \hbox{ in }\Omega\times (0,\,T)\\
\noalign{\smallskip}
y\vert_{\Sigma}=0\\
\noalign{\smallskip}
y(0)=y_0\hbox{ in }\Omega
\end{array}
\right.
$$
is in $L^{\infty}(0,\,T;L^4(\Omega)^N)$.
\end{theorem}
\textbf{Proof:}
To apply Theorem \ref{Hmenos1L4} to each component.
\begin{flushright}
$\square$
\end{flushright}

\begin{remark}
The results are true if we replace the operator $-\Delta$ by $-\nabla\cdot(A(t,x)\nabla)$ with $A\in L^\infty(\Omega\times (0,\,T))^{N\times N}$ and $\displaystyle{\sum_{i.j=1}^NA_{ij}(x,t)\xi_i\xi_j\ge \alpha\vert\xi\vert^2}$ for a. e. $(x,t)\in \Omega\times (0,\,T)$ and $\forall \xi\in \R^N$.
\end{remark}

\section{The uniqueness question.}
It is well-known that the Navier-Stokes equations are the classical model to study a newtonian fluid:
$$\left\{
\begin{array}{l}
y_t-\nu \Delta y+y\cdot \nabla y+\nabla p=f\hbox{ in }\Omega\times (0,\,T)\\
\noalign{\smallskip}
\nabla \cdot y=0\hbox{ in }\Omega\times (0,\,T)\\
\noalign{\smallskip}
y\vert_\Sigma=0\\
\noalign{\smallskip}
y(0)=y_0\hbox{ in }\Omega,
\end{array}
\right.
$$
being $\nu>0$ the cinematic viscosity, $y(x,t)$ the velocity of the fluid and $p(x,t)$ its pression.\newline
Following \cite{Teman}, we recall the usual spaces to treat these equations.\newline
Let be 
$${\cal{V}}:=\{ v\in {\cal{D}}(\Omega)^N:\ \ \nabla \cdot v=0\},$$
where ${\cal{D}}(\Omega)$ is the space of infinitely differentiable functions and with compact support in $\Omega$.\newline
We define $V$ as the closure of $\cal{V}$ in the topology of $H^1_0(\Omega)$, and it is characterized by
$$V:=\{v\in H^1_0(\Omega)^N:\ \ \nabla \cdot v=0\},$$
and $H$ is the closure of ${\cal{V}}$ in the topology of $L^2(\Omega)$,
$$H:=\{v\in L^2(\Omega)^N:\nabla \cdot v=0,\ \ v\cdot \vec{n}=0\hbox{ on }\partial \Omega\}$$
($\vec{n}$ normal vector, outside of $\partial \Omega$).\newline
In $H$ it is considered the scalar product induced by $L^2(\Omega)^N$ and in $V$ the scalar product is
$$((u,v)):=\sum_{i=1}^N(D_iu,D_iv),$$
being $(\cdot,\cdot)$ the scalar product in $L^2(\Omega)^N$.\newline
We denote $\vert \cdot\vert$ the norm in $L^2(\Omega)^N$ and $\Vert \cdot\Vert$ the usual norm in $H^1_0(\Omega)^N$.\newline
The space $V'$ is the dual of $V$ and its norm $\Vert \cdot \Vert_\ast$ is the associated to $\Vert\cdot\Vert$.
We have
$$V\subset H\subset V'.$$
where we have used the Riesz Theorem to identify $H$ and $H'$. The injections are compact (they would be only continuous if $\Omega$ were not bounded). Besides, they are dense. The duality product is denoted by $<\cdot,\cdot>$.
\vskip 0.2cm
A formulation of the problem is\newline
Find $y:\Omega\times (0,\,T)\to \R^N$ and $p:\Omega\times (0,\,T)\to \R$ such that
$$\left\{
\begin{array}{l}
y_t-\nu \Delta y+y\cdot \nabla y+\nabla p=f\hbox{ in }\Omega\times (0,\,T)\\
\noalign{\smallskip}
\nabla \cdot y=0\hbox{ in }\Omega\times (0,\,T)\\
\noalign{\smallskip}
y\vert_\Sigma=0\\
\noalign{\smallskip}
y(0)=y_0\hbox{ in }\Omega,
\end{array}
\right.
$$
with $f:\Omega\times (0,\,T)\to \R^N$ and $y_0:\Omega\to \R^N$ given.
\vskip 0.2cm
The definition of weak solution is given below.
\begin{definition}
Let be $f\in L^2(0,\,T;V')$ and $y_0\in H$. It is said that $y$ is a weak solution of the Navier-Stokes problem if 
\begin{equation}\label{Pb1}
\left\{
\begin{array}{l}
y\in L^2(0,\,T;V)\bigcap L^\infty(0,\,T;H),\\
\noalign{\smallskip}
\displaystyle{\frac{d}{dt}}(y,v)+\nu((y,v))+b(y,y,v)=<f,v>\ \ \ \forall v\in V\\
\noalign{\smallskip}
y(0)=y_0,
\end{array}
\right.
\end{equation}
where $b$ is the trilinear form on $V$ with values in $\R^N$ given by
$$b(u,v,w):=\displaystyle{\sum_{i=1}^N\int_{\Omega}u_iD_iv_jw_j}\, dx\ \ \ \forall u,v,w\in V.$$
\end{definition}
An equivalent formulation of (\ref{Pb1}) is the following:\newline
\begin{equation}\label{Pb2}
\left\{
\begin{array}{l}
y\in L^2(0,\,T;V)\bigcap L^\infty(0,\,T;H),\ y_t\in L^1(0,\,T;V')\\
\noalign{\smallskip}
y_t-\nu Ay+B(y,y)=f\\
\noalign{\smallskip}
y(0)=y_0,
\end{array}
\right.
\end{equation}
being the operators $A$ and $B$
$$A:V\to V'$$
$$<Au,v>=((u,v))\ \ \ \forall u,v\in V$$
$$B:V\to V'$$
$$<Bu,v>=b(u,u,v)\ \ \ \forall u,v\in V.$$
The existence of weak solution has been established as well as some uniqueness results (see \cite{Teman}, \cite{Lions}) :
\begin{theorem}
Let be $f\in L^2(0,\,T;V')$ and $y_0\in H$. Then, there is at least one weak solution $y$ of the Navier-Stokes problem. Besides, $y$ is weakly continuous in $[0,\,T]$ to $H$.
\end{theorem}
The issue of the uniqueness of solution has been solved partially. If the spatial dimension is two, $N=2$, it is known that there is a unique weak solution. However, when the dimension is three, $N=3$, the uniqueness is guaranteed in two situations:
\begin{itemize}
\item Either we consider the solutions in $L^8(0,\,T;L^4(\Omega))$ or more general, the solutions are in $L^s(0,\,T;L^r(\Omega))$ with $\frac{2}{s}+\frac{3}{r}\le 1$ if $\Omega$ is bounded and $\frac{2}{s}+\frac{3}{r} = 1$ if $\Omega$ is not bounded (see \cite{Lions} page 84).
\item Or $f\in L^\infty(0,\,T;H)$, $f'\in L^1(0,\,T;H)$, $y_0\in H^2(\Omega)\cap V$, and $\nu$ is large enough or $f$ and $y_0$ are small enough.
\end{itemize}

The stationary Stokes problem is the linearized stationary form of the Navier-Stokes equations:
\begin{equation}\label{Stokes}
\left\{
\begin{array}{l}
-\nu \Delta u+\nabla p=f \hbox{ in }\Omega\\
\noalign{\smallskip}
\nabla \cdot u=0\hbox{ in }\Omega\\
\noalign{\smallskip}
u=0\hbox{ on }\partial \Omega,
\end{array}
\right.
\end{equation}
where $f\in H^{-1}(\Omega)^N$.\newline
A function $u\in H^1_0(\Omega)^N$ satisfying (\ref{Stokes}) in the following sense: there exists $p\in L^2(\Omega)$ such that
$$-\nu \Delta u+\nabla p=f$$
in the distribution sense in $\Omega$ and 
$$\nabla \cdot u=0$$
in the distribution sense in $\Omega$, is called a weak solution of the stationary Stokes problem (\ref{Stokes}).\newline
The equivalent variational formulation is:
$$u\in V\hbox{ such that }$$
$$\nu((u,v))=<f,v>\ \ \ \forall v\in V.$$
This problem has a unique solution. \newline
The mapping
$$\Lambda:f\in L^2(\Omega)^N \mapsto 1/\nu u,$$
being $u$ the solution of the stationary Stokes problem, is a compact linear operator and also self-adjoint, so the eigenfunctions of the operator constitute an orthonormal basis in $L^2(\Omega)^N$:
$$w_j\in V\ \ \ \ ((w_j,v))=\lambda_j(w_j,v)\ \ \ \forall v\in V.$$
\vskip 0.2cm

Our reasoning shall prove that any weak solution of the Navier-Stokes problem ( with $f\in L^2(0,\,T;H^{-1}(\Omega)^N)$ and $y_0\in H\cap L^4(\Omega)^N$) is in $L^\infty(0,\,T;L^4(\Omega)^N)$, and so, this solution is unique. After that, by a discretization procedure, we obtain the uniqueness for $y_0\in H$ and $f\in L^2(0,\,T;H^{-1}(\Omega)^N)$. Finally, when $f$ is more ge\-ne\-ral, $f\in L^2(0,\,T;V')$, we extend this operator to $L^2(0,\,T;H^{-1}(\Omega)^N)$ where the uniqueness has been already proved.
\begin{theorem}\label{unicidadeHmenos1L4}
The Navier-Stokes problem in three dimension, with $f\in L^2(0,\,T;H^{-1}(\Omega)^N)$ and $y_0\in H\cap L^4(\Omega)^N$, has a unique solution.
\end{theorem}
\textbf{Proof:}
We call $z$ the solution of 
$$\left\{
\begin{array}{l}
z_t-\nu\Delta z=0\hbox{ in }\Omega\times (0,\,T)\\
\noalign{\smallskip}
z\vert_{\Sigma}=0\\
\noalign{\smallskip}
z(0)=y_0/2\hbox{ in }\Omega,
\end{array}
\right.
$$
and let be 
$$w=y-z, $$
with $y$ a solution of Navier-Stokes problem.\newline
Then, $w$ is the solution of a linear parabolic problem:
$$
\left\{
\begin{array}{l}
w_t-\nu \Delta w=f-B(y,y)\\
\noalign{\smallskip}
w\vert_\Sigma=0\\
\noalign{\smallskip}
w(0)=y_0/2.
\end{array}
\right.
$$
The operator $B(y,y)$ is defined on $H^{-1}(\Omega)^N$ like this
$$<B(y,y),v>=b(y,y,v).$$
Then, $w$ is the solution of a linear system satisfying the conditions of Theorem \ref{sistemaHmenos1L4}, so $w$ is in $L^\infty(0,\,T;L^4(\Omega)^N)$. Since $y=w+z$,  $y\in L^\infty(0,\,T;L^4(\Omega)^N)$ and then, the Navier-Stokes solution has the sufficient regularity to guarantee it is unique (see \cite{Teman}).
\begin{flushright}
$\square$
\end{flushright}

\vskip 0.3cm
Next, we study the uniqueness when the initial data is in $L^2(\Omega)^N$ and $f$ is still in $L^2(0,\,T;H^{-1}(\Omega)^N)$ .
\begin{theorem}\label{unicidadenL2}
The Navier-Stokes problem in three dimension, with $f\in L^2(0,\,T;H^{-1}(\Omega)^N)$ and $y_0\in H$ has a unique solution.
\end{theorem}
\textbf{Proof:}
Let be $\{ w_j\}_j$ the sequence of eigenfunctions of the Stokes problem.
We define
$$(y_0^{(n)})_k=\sum_{j=1}^n(y_{0k},w_j)w_j,$$
being $y_{0k}$ the k-th component of $y_0$.\newline
We call $y^{(n)}$ the solution of the Navier-Stokes problem with initial data ${y_0}^{(n)}$:
\begin{equation}\label{NSn}
\left\{
\begin{array}{l}
\displaystyle{\frac{d}{dt}}(y^{(n)},v)+\nu((y^{(n)},v))+b(y^{(n)},y^{(n)},v)=<f,v> \ \ \ v\in V\\
\noalign{\smallskip}
y^{(n)}(0)={y_0}^{(n)}.
\end{array}
\right.
\end{equation}
Since the initial data is in $L^4(\Omega)^N$ and the right-hand side is in $L^2(0,\,T;H^{-1}(\Omega)^N)$, by Theorem \ref{unicidadeHmenos1L4}, this solution exists and it is unique. Besides, it is in $L^\infty(0,\,T;L^4(\Omega))$.\newline
Let be $y$ a solution of Navier-Stokes problem with initial data $y_0$. We shall prove that the sequence $\{y^{(n)}\}$ converges to $y$. Since $\{y^{(n)}\}_n$ is uniquely determined there cannot be two different solutions of Navier-Stokes problem.
\newline
Subtracting the equations satisfied by $y$ and $y^{(n)}$ we obtain
$$\displaystyle{\frac{d}{dt}}(y-y^{(n)},v)+\nu ((y-y^{(n)},v))+b(y,y,v)-b(y^{(n)},y^{(n)},v)=0\ \ \ \forall v\in V.$$
Taking $y(t)-y^{(n)}(t)$ as a test function and integrating in $\Omega$ we obtain
$$\displaystyle{\frac{d}{dt}}\vert y(t)-y^{(n)}(t)\vert^2+\nu\Vert y(t)-y^{(n)}(t)\Vert ^2+b(y,y,y-y^{(n)})-b(y^{(n)},y^{(n)},y-y^{(n)})=0.$$
Using the trilinearity and the properties $b(u,v,v)=0$ and $b(u,v,w)=-b(u,w,v)$ we have
$$b(y,y,y-y^{(n)})-b(y^{(n)},y^{(n)},y-y^{(n)})=-b(y,y,y^{(n)})-b(y^{(n)},y^{(n)},y)=b(y-y^{(n)},y^{(n)},y).$$
Then, 
\begin{equation}\label{igualdadb1}
\displaystyle{\frac{d}{dt}}\vert y(t)-y^{(n)}(t)\vert^2+\nu\Vert y(t)-y^{(n)}(t)\Vert ^2+b(y-y^{(n)},y^{(n)},y)=0.
\end{equation}
Again, by the properties of $b$,
$$b(y-y^{(n)},y^{(n)},y)=b(y-y^{(n)},y^{(n)}-y,y)=-b(y-y^{(n)},y-y^{(n)},y).$$
Applying this equality to (\ref{igualdadb1}) and the Holder inequality, we have
$$\displaystyle{\frac{d}{dt}}\vert y(t)-y^{(n)}(t)\vert^2+\nu\Vert y(t)-y^{(n)}(t)\Vert ^2=b(y-y^{(n)},y-y^{(n)},y)\le$$
$$\le \int_{\Omega}\vert y-y^{(n)}\vert \,\vert \nabla (y-y^{(n)})\vert\, \vert y\vert\le \Vert y(t)-y^{(n)}(t)\Vert _{L^4(\Omega)}\Vert \nabla(y(t)-y^{(n)}(t))\Vert_{L^2(\Omega)}\Vert y(t)\Vert_{L^4(\Omega)}.$$
When the spatial dimension is three there is an estimation of the norm in $L^4(\Omega)$, for any open set $\Omega$ and for any $v\in H^1_0(\Omega)$ (see Lemma 3.5, page 296 in \cite{Teman}). It is hold
$$\Vert v\Vert_{L^4(\Omega)}\le 2^{1/2}\Vert v\Vert_{L^2(\Omega)}^{1/4}\Vert \nabla v\Vert_{L^2(\Omega)}^{3/4}.$$
We apply this inequality to $y(t)-y^{(n)}(t)$. Then,
$$\Vert y(t)-y^{(n)}(t)\Vert_{L^4(\Omega)}\Vert\nabla (y(t)-y^{(n)}(t))\Vert_{L^2(\Omega)}\Vert y(t)\Vert_{L^4(\Omega)}\le $$
$$\le 2^{1/2}\Vert y(t)-y^{(n)}(t)\Vert_{L^2(\Omega)}^{1/4}\Vert\nabla (y(t)-y^{(n)}(t))\Vert_{L^2(\Omega)}^{3/4}\Vert\nabla (y(t)-y^{(n)}(t))\Vert_{L^2(\Omega)}\Vert y(t)\Vert_{L^4(\Omega)}=$$
$$=2^{1/2}\Vert y(t)\Vert_{L^4(\Omega)}\Vert y(t)-y^{(n)}(t)\Vert_{L^2(\Omega)}^{1/4}\Vert\nabla (y(t)-y^{(n)}(t))\Vert_{L^2(\Omega)}^{7/4}.$$
By Young's inequality with $p=8$ and $q=8/7$ we have
$$2^{1/2}\Vert y(t)\Vert_{L^4(\Omega)}\Vert y(t)-y^{(n)}(t)\Vert_{L^2(\Omega)}^{1/4}\Vert\nabla (y(t)-y^{(n)}(t))\Vert_{L^2(\Omega)}^{7/4}\le$$
$$\le C\displaystyle{(\Vert y(t)-y^{(n)}(t)\Vert_{L^2(\Omega)}^{1/4})^8+\frac{\nu}{2}(\Vert\nabla(y(t)-y^{(n)}(t))\Vert_{L^2(\Omega)}^{7/4})^{8/7}}=$$
$$=C\displaystyle{\Vert y(t)-y^{(n)}(t)\Vert^2_{L^2(\Omega)}+\frac{\nu}{2}\Vert\nabla(y(t)-y^{(n)}(t))\Vert^2_{L^2(\Omega)}}.$$
Then,
$$\displaystyle{\frac{d}{dt}}\vert y(t)-y^{(n)}(t)\vert^2+\nu\Vert y(t)-y^{(n)}(t)\Vert ^2\le C\displaystyle{\vert y(t)-y^{(n)}(t)\vert^2+\frac{\nu}{2}\Vert\nabla(y(t)-y^{(n)}(t))\Vert^2}.$$
So,
$$\displaystyle{\frac{1}{2}\frac{d}{dt}\vert y(t)-y^{(n)}(t)\vert^2}\le C\vert y(t)-y^{(n)}(t)\vert^2.$$
By Gronwall's Lemma,
$$\vert y(t)-y^{(n)}(t)\vert^2\le e^{Ct}\vert y_0-{y_0}^{(n)}\vert^2,$$
Since $\{y_0^{(n)}\}$ converges to $y_0$ in $L^2(\Omega)^N$ we have
$$\vert y(t)-y^{(n)}(t)\vert\to 0,$$
and this yields that
$$y^{(n)}\to y\hbox{ in }L^\infty(0,\,T;L^2(\Omega)).$$ 
This proves that the solution of Navier-Stokes problem is unique since for any solution there is the same sequence of functions that converges to it.
\begin{flushright}
$\square$
\end{flushright}

Now, let us relax the hypothesis about $f$. 
\begin{theorem}\label{unicidadVprima}
The Navier-Stokes problem in three dimension, with $f\in L^2(0,\,T;V')$ and $y_0\in H$ has a unique solution.
\end{theorem}
\textbf{Proof:}
Let be $\tilde{f}\in L^2(0,\,T;H^{-1}(\Omega)^N)$ given by
$$<\tilde{f}(t),v>=<f(t),\overline{v}>+<q,v-\overline{v}>\ \ \ \forall v\in H^1_0(\Omega)^N,$$
being $\overline{v}$ the projection of $v$ in $V$ and $q$ any function in $L^2(\Omega)^N$.\newline
Then, $\tilde{f}$ is in $L^2(0,\,T;H^{-1}(\Omega)^N$. Applying Theorem \ref{unicidadenL2}, the Navier-Stokes problem
$$
\left\{
\begin{array}{l}
y_t-\nu Ay+B(y,y)=\tilde{f}\\
\noalign{\smallskip}
y(0)=y_0,
\end{array}
\right.
$$
has a unique solution. This means that there exists a unique $y$ verifying
$$
\left\{
\begin{array}{l}
y\in L^2(0,\,T;V)\bigcap L^\infty(0,\,T;H),\\
\noalign{\smallskip}
\displaystyle{\frac{d}{dt}}(y,v)+\nu((y,v))+b(y,y,v)=<\tilde{f},v>\ \ \ \forall v\in V\\
\noalign{\smallskip}
y(0)=y_0.
\end{array}
\right.
$$
But, when $v\in V$, $\tilde{f}=f$ so, for $f\in L^2(0,\,T;V')$ and $y_0\in H$ there exists a unique weak solution.
\begin{flushright}
$\square$
\end{flushright}

\section{Conclusions}

In this work, we have obtained several significant results. First, we have established a connection between the solution of the initial boundary value problem for the heat equation and the solution of any linear parabolic problem. This relation provides a useful framework to analyze the behavior of solutions to more general parabolic equations. Second, we have proved a regularity result for the solutions of linear parabolic problems: if the initial data belongs to $L^4(\Omega)$, then the solution remains in $L^\infty(0,,T;L^4(\Omega))$. This gain in regularity plays a crucial role in the final part of the work, where we address the three-dimensional Navier-Stokes equations. By relying on the tools developed throughout the paper, we have been able to prove the uniqueness of the solution to the Navier-Stokes problem in three dimensions under the same initial regularity assumption.

As for future research directions, it would be interesting to analyze whether similar techniques can be extended to other classes of nonlinear parabolic equations or to problems with more general boundary conditions. Another possible line of work is to explore the numerical implementation of the regularity result obtained for linear parabolic equations and its implications in the simulation of fluid dynamics. Finally, one could investigate whether the approach presented here can be adapted to treat more general fluid-structure interaction problems.  






\vskip 12cm
\begin{center}
\scriptsize{\lq\lq You have made us for Yourself, O Lord, and our heart is restless till it rests in You\rq\rq.  Saint Augustine.}
\end{center}

\end{document}